\documentclass[12pt]{article}

\usepackage{amsmath, amssymb} 
\usepackage{subcaption}
\usepackage{amsthm} 
\usepackage{enumitem}  
\usepackage{tikz,tkz-graph}
\usepackage{float}
\usepackage{graphicx}
\usepackage{cite}
\usepackage{blkarray}
\usepackage{authblk}
\usepackage{hyperref}

\DeclareMathOperator\rk{rk}

\DeclareMathOperator\col{Col}

\DeclareMathOperator\SPAN{span}

\newtheoremstyle{plainsl}%
	{\topsep}
	{\topsep}
	{\slshape} 
	{}
	{\normalfont\bfseries}
	{.}
	{ }
	{}

\swapnumbers

{\theoremstyle{plainsl}
\newtheorem{theorem}{Theorem}[section]

\newtheorem{lemma}[theorem]{Lemma}
\newtheorem{corollary}[theorem]{Corollary}}
{\theoremstyle{remark}
}

\newcommand\cref[1]{Corollary~\ref{cor:#1}}

\renewcommand\proof{\noindent\textsl{Proof. }}
\newcommand\sqr[2]{{\vbox{\hrule height.#2pt
    \hbox{\vrule width.#2pt height#1pt \kern#1pt
        \vrule width.#2pt}\hrule height.#2pt}}}

\renewcommand\qed{%
	\ifmmode\eqno\sqr53
	\else\nolinebreak\ \hfill\sqr53\medbreak\fi}

\numberwithin{equation}{section}

\newcommand{\widehatP}{\widehat{P}}
\newcommand{\abs}[1]{\left| #1 \right|} 
\newcommand{\Mod}[1]{\ (\mathrm{mod}\ #1)}

\begin{document}
\title{Hamiltonians of Bipartite Walks}

\author[1]{Qiuting Chen}
\author[1]{Chris Godsil}
\author[1]{Mariia Sobchuk}
\author[2]{Hanmeng Zhan}
\affil[1]{Department of Combinatorics \& Optimization,
	University of Waterloo, Waterloo, Ontario, Canada}
\affil[2]{Department of Mathematics and Statistics, York University, Toronto, Ontario, Canada}

\maketitle

\pagenumbering{arabic}

\maketitle

\abstract{In this paper, we introduce a discrete quantum walk model called bipartite walks. Bipartite walks include many known discrete quantum walk models, like arc-reversal walks, vertex-face walks. For the transition matrix of a quantum walk, there is a Hamiltonian associated with it. We will study the Hamiltonians of the bipartite walks. Let $S$ be a skew-symmetric matrix. We are mainly interested in the Hamiltonians  of the form $iS$. We show that the Hamiltonian can be written as $iS$ if and only if the adjacency matrix of the bipartite graph is invertible. We show that arc-reversal walks and vertex-face walks are special cases of bipartite walks. Via the Hamiltonians, phenomena of bipartite walks lead to phenomena of continuous walks. We show in detail how we use bipartite walks on paths to construct universal perfect state transfer in continuous walks. }

\section{Introduction}

Quantum walks are a quantum mechanical analogue of classical random walks. They provide a powerful tool for the study and development of quantum algorithms~\cite{Childs2002,szegedy2004}. Based on how time evolves, a quantum walk can be either continuous or discrete. For discrete quantum walks, there are several models that have been proposed and studied~\cite{Staggered,coinedwalk,szegedy2004}. In this paper, the walks we focus on are called  bipartite walks; they generalize many known models such as arc-reversal walks and vertex-face walks. 

We turn to a description of bipartite walks. A discrete quantum walk is given by a unitary operator
$U$ on a complex vector space $\mathbb{C}^n$. 
We refer to $U$ as the \textsl{transition matrix} of a discrete quantum walk.
The state of the underlying quantum system is a unit
vector in $\mathbb{C}^n$. If the initial state is $z$, then after $k$ steps of the walk, the state
is $U^kz$. This is a unit vector, and so the squared absolute values of its entries sum to $1$.
The outcome of a measurement after $k$ steps is an element $i$ of $\{1,\ldots,n\}$, and the
probabilty that the result is $i$ is $|(U^kz)_i|^2$.

In our case, the state space is the space of complex functions
on the edges of a bipartite graph $G$. We assume that $X$ and $Y$ are the two colour classes 
of $G$ and using these we construct two partitions of $E(G)$.
For the first partition, $\pi_0$, two edges are in the 
same cell if they have a vertex in common, and that vertex is in $X$. 
For the second partition
$\pi_1$, two edges are in the same cell if they have a vertex in common, and that vertex is in $Y$. Each of these partitions determines a projection, namely the projection onto the
functions on $E(G)$ that are constant on the cells of $\pi_0$ and $\pi_1$.
We denote these projections by $P$ and $Q$ respectively.

If $R$ is a projection, then
\[
	(2R-I)^2 = 4R^2 -4R +I =4R -4R +I = I
\]
and, since $R=R^*$, we see that $2R-I$ is unitary. (Geometrically it is a reflection.) 
Hence we can define a unitary operator $U$ by
\[
	U := (2P-I)(2Q-I).
\]
This the transition matrix of the bipartite walk on $G$.

Konno et al.~in~\cite{twopartition}  introduce a family of discrete-time quantum walks, called two-partition model, which is based on two equivalence-class partitions of the computational basis. The two partition used in the two-partition model does not necessarily give us two reflections.
Bipartite walks are a special case of the two-partition model introduced by 
Konno et al.~in~\cite{twopartition}. Note that the paper by Konno et al.~focuses on showing the unitary equivalence between the members of two-partition model while we study the Hamiltonian of the transition matrix of the bipartite walk in this paper.

On the other hand, many of the most
commonly used discrete walks can be formulated as bipartite walks. We will give a constructive proof to show that arc-reversal walk can be viewed as a special case of bipartite walk.

There is a second class of quantum walks: \textsl{continuous quantum walks}. Here
the state space is the space of complex functions on the vertices of a graph $G$.
The walk is specified by a Hermitian matrix $H$ with rows and columns indexed by the
vertices of $G$ (for example, the adjcency matrix of $G$). We then
define transition matrices $U(t)$ by
\[
	U(t) := \exp(itH),\quad(t \in \mathbb{R}).
\]
If the initial state of the walk is given by the unit vector $z$, the state at time $t$
is $U(t)z$. 
For each unitary matrix $U$, there are Hermitian matrices $H$ such that 
\[
	U=\exp(iH).
\] 
(We refer to $H$ as a \textsl{Hamiltonian} of $U$.) 
It follows that a discrete walk on $G$ gives rise to a continuous quantum walk on the
edges of $G$ and if the continuous walk is given by matrices $U(t)$, the transition
matrix for the discrete walk is $U(1)$.

Our goal in this paper is to study the Hamiltonians of bipartite wallks. This is a topic that
has not been studied before.

For the discrete quantum walk governed by the unitary matrix $U$, there is a Hamiltonian $H$ associated with it. When there is a real skew-symmetric $S$ such that the Hamiltonian $H$ is of $H=iS$, it can be viewed as the skew-adjacency matrix of a oriented weighted graph, which we call \textsl{the $H$-digraph}. Hamiltonians of quantum walks are often associated with continuous quantum walks and have not been considered in the context of discrete quantum walks.

 So far, most studies of the bipartite walk have been limited to the transition matrix and the behaviors of the walk~\cite{szegedy2004,twopartition,StefanakSkoupy}.
In this paper, we study Hamiltonians of bipartite walks and $H$-digraphs associated with it. 
Spectral properties of the transition matrix is the main tool we exploit to study the Hamiltonian of $U$. 

Let $S$ be a skew-symmetric matrix. We are mainly interested in the case when the Hamiltonian $H$ can be written as $H=iS$, which is not always true.We prove that the Hamiltonian $H$ is of the form $H=iS$ if and only if the adjacency matrix of $G$ is invertible. 

As mentioned before, vertex-face walk can be viewed as a special case of bipartite walk. In Section~\ref{vf walk}, we show the equivalence relations between bipartite walks and vertex-face walks. The Hamiltonians obtained from vertex-face walks have some interesting properties, which have been studied extensively in~\cite{harmonyphd}. Here we introduce those properties and rephrase them from perspective of bipartite walk in Section~\ref{vf walk} and Section~\ref{vxf on CompleteG}.

When $G$ is a path on $n$ vertices, the transition matrix of the bipartite walk is a permutation matrix. When $n\geq 4$ is even, the associated $H$-digraph is a weighted oriented $K_{n-1}$. When $n\equiv 3\Mod 4$, the associated $H$-digraph is two copies of a weighted oriented $K_\frac{n-1}{2}$. Similar results can also be proved for the bipartite walk on even cycles. 

Studying the Hamiltonian of bipartite walks helps us to construct examples of continuous walks with desired properties. Consider continuous quantum walk on a graph $G$ and the Hamiltonian is the adjacency matrix of $G$. If the walk has perfect state transfer between every pair of vertices of $G$, the walk has universal perfect state transfer. This is a rare and interesting phenomenon.
Using the properties of bipartite walks on paths and cycles, we find a way to weight the edges of complete graphs such that the resulting weighted graph has universal perfect state transfer. This demonstrates how we can use the Hamiltonian and bipartite walks to construct some interesting but previously hard-to-find phenomenon in continuous walks.

\section{Preliminaries}
\label{Intro}

Let $G$ be a $(d_0,d_1)$-biregular bipartite graph with two parts $C_0,C_1$. 
Now we define two partitions of the edges of $G$, denoted by $\pi_0,\pi_1$ respectively.  If two edges have the same end $x$ in $C_0$, then they belong to the same cell of $\pi_0$. Similarly, if two edges have the same end $y$ in $C_1$, then they belong to the same cell of $\pi_1$.

Given a matrix $M$, we \textsl{normalize} it by scaling each column of $M$ to a unit vector.
Let $P_0,P_1$ be characteristic matrix of $\pi_0,\pi_1$ respectively and let $\widehatP_0,\widehatP_1$ denote the normalized $P_0,P_1$ respectively.

Let 
\[
	P=\widehatP_0\widehatP_0^T,\quad Q=\widehatP_1\widehatP_1^T
\] 
be the projections onto the vectors that is constant on the cells of $\pi_0,\pi_1$ respectively. We define the transition matrix of the bipartite walk over $G$ to be
\[
	U=\left(2\widehatP_0\widehatP_0^T-I\right)\left(2\widehatP_1\widehatP_1^T-I\right)
		=\left(2P-I\right)\left(2Q-I\right).
\] 

\begin{figure}[H]
\begin{center}

\begin{tikzpicture}
\definecolor{cv0}{rgb}{0.0,0.0,0.0}
\definecolor{cfv0}{rgb}{1.0,1.0,1.0}
\definecolor{clv0}{rgb}{0.0,0.0,0.0}
\definecolor{cv1}{rgb}{0.0,0.0,0.0}
\definecolor{cfv1}{rgb}{1.0,1.0,1.0}
\definecolor{clv1}{rgb}{0.0,0.0,0.0}
\definecolor{cv2}{rgb}{0.0,0.0,0.0}
\definecolor{cfv2}{rgb}{1.0,1.0,1.0}
\definecolor{clv2}{rgb}{0.0,0.0,0.0}
\definecolor{cv3}{rgb}{0.0,0.0,0.0}
\definecolor{cfv3}{rgb}{1.0,1.0,1.0}
\definecolor{clv3}{rgb}{0.0,0.0,0.0}
\definecolor{cv4}{rgb}{0.0,0.0,0.0}
\definecolor{cfv4}{rgb}{1.0,1.0,1.0}
\definecolor{clv4}{rgb}{0.0,0.0,0.0}
\definecolor{cv5}{rgb}{0.0,0.0,0.0}
\definecolor{cfv5}{rgb}{1.0,1.0,1.0}
\definecolor{clv5}{rgb}{0.0,0.0,0.0}
\definecolor{cv6}{rgb}{0.0,0.0,0.0}
\definecolor{cfv6}{rgb}{1.0,1.0,1.0}
\definecolor{clv6}{rgb}{0.0,0.0,0.0}
\definecolor{cv7}{rgb}{0.0,0.0,0.0}
\definecolor{cfv7}{rgb}{1.0,1.0,1.0}
\definecolor{clv7}{rgb}{0.0,0.0,0.0}
\definecolor{cv0v1}{rgb}{0.0,0.0,0.0}
\definecolor{cv0v5}{rgb}{0.0,0.0,0.0}
\definecolor{cv1v2}{rgb}{0.0,0.0,0.0}
\definecolor{cv1v4}{rgb}{0.0,0.0,0.0}
\definecolor{cv2v3}{rgb}{0.0,0.0,0.0}
\definecolor{cv5v6}{rgb}{0.0,0.0,0.0}
\definecolor{cv6v7}{rgb}{0.0,0.0,0.0}
\Vertex[style={minimum size=1.0cm,draw=cv0,fill=cfv0,text=clv0,shape=circle},LabelOut=false,L=\hbox{$0$},x=0cm,y=4cm]{v0}
\Vertex[style={minimum size=1.0cm,draw=cv1,fill=cfv1,text=clv1,shape=circle},LabelOut=false,L=\hbox{$1$},x=4cm,y=4cm]{v1}
\Vertex[style={minimum size=1.0cm,draw=cv2,fill=cfv2,text=clv2,shape=circle},LabelOut=false,L=\hbox{$2$},x=0cm,y=3cm]{v2}
\Vertex[style={minimum size=1.0cm,draw=cv3,fill=cfv3,text=clv3,shape=circle},LabelOut=false,L=\hbox{$3$},x=4cm,y=3cm]{v3}
\Vertex[style={minimum size=1.0cm,draw=cv4,fill=cfv4,text=clv4,shape=circle},LabelOut=false,L=\hbox{$4$},x=0.0cm,y=2cm]{v4}
\Vertex[style={minimum size=1.0cm,draw=cv5,fill=cfv5,text=clv5,shape=circle},LabelOut=false,L=\hbox{$5$},x=4cm,y=2cm]{v5}
\Vertex[style={minimum size=1.0cm,draw=cv6,fill=cfv6,text=clv6,shape=circle},LabelOut=false,L=\hbox{$6$},x=0cm,y=1cm]{v6}
\Vertex[style={minimum size=1.0cm,draw=cv7,fill=cfv7,text=clv7,shape=circle},LabelOut=false,L=\hbox{$7$},x=4cm,y=1.0cm]{v7}
\Edge[lw=0.1cm,style={color=cv0v1,},](v0)(v1)
\Edge[lw=0.1cm,style={color=cv0v5,},](v0)(v5)
\Edge[lw=0.1cm,style={color=cv1v2,},](v1)(v2)
\Edge[lw=0.1cm,style={color=cv1v4,},](v1)(v4)
\Edge[lw=0.1cm,style={color=cv2v3,},](v2)(v3)
\Edge[lw=0.1cm,style={color=cv5v6,},](v5)(v6)
\Edge[lw=0.1cm,style={color=cv6v7,},](v6)(v7)
\end{tikzpicture}
\caption{Bipartite graph on $8$ vertices}
\label{not return pst graph}
\end{center}
\end{figure}
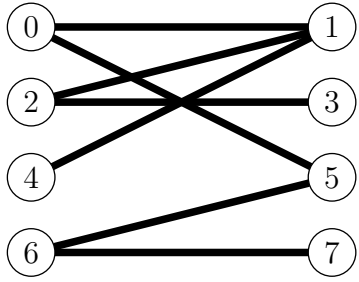

Now consider the bipartite graph $G$ in Figure~\ref{not return pst graph} as an example. We define a bipartite walk on $G$. The two parts of $G$ are $C_0=\{0,2,4,6\}$ and $C_2=\{1,3,5,7\}$. For the partitions $\pi_0,\pi_1$, the edge $(0,1),(0,5)$ are in the same cell in $\pi_0$ and Edge $(0,1),(2,1),(4,1)$ are in the same cell in $\pi_1$.
We have that  
\[
	\hat{P}_0=\begin{pmatrix}
\frac{1}{\sqrt{3}} & 0 & 0 & 0 \\
0 & 0 & \frac{1}{\sqrt{2}} & 0 \\
\frac{1}{\sqrt{3}} & 0 & 0 & 0 \\
\frac{1}{\sqrt{3}}& 0 & 0 & 0 \\
0 & 1 & 0 & 0 \\
0 & 0 & \frac{1}{\sqrt{2}}& 0 \\
0 & 0 & 0 & 1
\end{pmatrix} ,\quad
\hat{P}_1=\begin{pmatrix}
 \frac{1}{\sqrt{2}} & 0 & 0 & 0 \\
 \frac{1}{\sqrt{2}}& 0 & 0 & 0 \\
0 & \frac{1}{\sqrt{2}}& 0 & 0 \\
0 & 0 & 1 & 0 \\
0 &  \frac{1}{\sqrt{2}}& 0 & 0 \\
0 & 0 & 0 &  \frac{1}{\sqrt{2}}\\
0 & 0 & 0 &  \frac{1}{\sqrt{2}}
\end{pmatrix}
\]
and hence, the corresponding projections are
\[
P=
\begin{pmatrix}
\frac{1}{3}& 0 & \frac{1}{3}& \frac{1}{3}& 0 & 0 & 0 \\[2.5mm]
0 & \frac{1}{2} & 0 & 0 & 0 & \frac{1}{2}& 0 \\[2.5mm]
\frac{1}{3}& 0 & \frac{1}{3} & \frac{1}{3}& 0 & 0 & 0 \\[2.5mm]
\frac{1}{3}& 0 & \frac{1}{3}&\frac{1}{3} & 0 & 0 & 0 \\[2.5mm]
0 & 0 & 0 & 0 & 1 & 0 & 0 \\[2.5mm]
0 & \frac{1}{2}& 0 & 0 & 0 & \frac{1}{2} & 0 \\[2.5mm]
0 & 0 & 0 & 0 & 0 & 0 & 1
\end{pmatrix},\quad
Q=\begin{pmatrix}
\frac{1}{2}& \frac{1}{2} & 0 & 0 & 0 & 0 & 0 \\[2.5mm]
\frac{1}{2}& \frac{1}{2}& 0 & 0 & 0 & 0 & 0 \\[2.5mm]
0 & 0 & \frac{1}{2}& 0 & \frac{1}{2} & 0 & 0 \\[2.5mm]
0 & 0 & 0 & 1 & 0 & 0 & 0 \\[2.5mm]
0 & 0 & \frac{1}{2} & 0 & \frac{1}{2} & 0 & 0 \\[2.5mm]
0 & 0 & 0 & 0 & 0 & \frac{1}{2} & \frac{1}{2} \\[2.5mm]
0 & 0 & 0 & 0 & 0 & \frac{1}{2} & \frac{1}{2}
\end{pmatrix}.
\]
The transition matrix of the bipartite walk on $G$ is 
\[
	U=\begin{pmatrix}
0 & -\frac{1}{3} & 0 &  \frac{2}{3}&  \frac{2}{3} & 0 & 0 \\[2.5mm]
0 & 0 & 0 & 0 & 0 & 0 & 1 \\[2.5mm]
0 & \frac{2}{3}  & 0 & \frac{2}{3}  & \frac{1}{3} & 0 & 0 \\[2.5mm]
0 &  \frac{2}{3}  & 0 & -\frac{1}{3}  &  \frac{2}{3}  & 0 & 0 \\[2.5mm]
0 & 0 & 1 & 0 & 0 & 0 & 0 \\[2.5mm]
1 & 0 & 0 & 0 & 0 & 0 & 0 \\[2.5mm]
0 & 0 & 0 & 0 & 0 & 1 & 0
\end{pmatrix}.
\]

Let $C$ denote the characteristic matrix of the incidence relation between $\pi_0,\pi_1$ with its rows indexed by the cells of $\pi_1$ and its columns indexed by the cells of $\pi_0$ such that 
\[
	C_{i,j}=1
\] 
if there is an edge that belongs to both $c_i$ in $\pi_1$ and $c_j$ in $\pi_0$. Then we have that 
\[
	C=P_1^TP_0
\] 
and normalized $C$ is 
\[
	\hat{C}=\widehatP_1^T\widehatP_0.
\]
The adjacency matrix of $G$ can be written as
\[
\quad A(G)=
\begin{pmatrix}
\mathbf{0}&C\\
C^T&\mathbf{0}
\end{pmatrix}.
\]

The incidence matrix and the normalized incidence matrix of the bipartite graph in Figure~\ref{not return pst graph} are 
\[
C=
\begin{pmatrix}
1 & 0 & 1 & 0 \\
1 & 1 & 0 & 0 \\
1 & 0 & 0 & 0 \\
0 & 0 & 1 & 1
\end{pmatrix},\quad
\hat{C}=\begin{pmatrix}
\frac{1}{\sqrt{6}} & 0 & \frac{1}{2} & 0 \\[3mm]
\frac{1}{\sqrt{6}} &\frac{1}{\sqrt{2}} & 0 & 0 \\[3mm]
\frac{1}{\sqrt{3}}& 0 & 0 & 0 \\[3mm]
0 & 0 & \frac{1}{2}& \frac{1}{\sqrt{2}}
\end{pmatrix}.
\]

\section{Arc-reversal walks are a special case}

Arc-reversal walks are a well-studied model and in this section, we give a constructive proof that arc-reversal walks can be considered as a special case of bipartite walks.

Given a graph $G$, we show that the bipartite walk on the subdivision graph of $G$ is equivalent to the arc-reversal walk on $G$.

For a graph $G$, we define a new graph $G'$ by subdivided every edge of $G$ and we call $G'$ \textsl{the subdivision graph of $G$}. Then $G'$ is a bipartite graph with parts $C_0=V(G')\backslash V(G)$ and $C_1=V(G)$. We define a bipartite walk on $G'$ with transition matrix 
\[
	U=(2P-I)(2Q-I).
\] For each vertex $\alpha\in C_0$ and $a\in C_1$, we have
 \[
	\deg_{G'}(\alpha)=2, \quad\deg_{G'}(a)=\deg_{G}(a).
\]

Now if  every edge $e$ of $G$ is replaced by two arcs $e_1,e_2$ with opposite directions, we can view the subdivision graph $G'$ as directed graph of $G$. Every edge in $G'$ can be viewed as an arc of directed $G$. 

Let
 \[
 	G_a=\frac{1}{\deg(a)}J-I
 	\]
be the Grover coin associated with vertex $a$. Then we have that 
\[
	2Q-I=\bigoplus_{v\in C_1} G_v=\begin{pmatrix}
	G_{v_1}&&&\\
	&G_{v_2}&&\\
	&&\ddots&\\
	&&&G_{v_n}
	\end{pmatrix}, 
	\]
where we assign the Grover coin to $v_i$ for every vertex $v_i$ in $V(G)$. Also, we have that
 \[
 	2P-I=\bigoplus_{v\in C_0}\frac{1}{2}J_2-I=\begin{pmatrix}
	\frac{1}{2}J_2-I&&&\\[3mm]
	&\frac{1}{2}J_2-I&&\\[3mm]
	&&\ddots&\\[3mm]
	&&&\frac{1}{2}J_2-I
	\end{pmatrix},
 \]
 which can be viewed as the arc-reversal matrix $R$, i.e., 
 \[
	R\cdot (a,b)=(b,a)
\]
for every arc $(a,b)$. Thus, every bipartite walk defined on the subdivision graph of $G$ is equivalent to the arc-reversal walk on $G$.

\section{Spectrum of transition matrix $U$}

Spectral properties of the transition matrix $U$ are the main machinery that we use to analyse the Hamiltonian of $U$. In this section, we present a complete characterization on the eigenvalues and eigenspaces of $U$. All the statements presented here are proved in \cite{harmonyphd} by Zhan in detail, so in this paper we omit the proofs. Note that here we use the same notations as defined before and so,
\[
	P=\widehatP_0\widehatP_0^T,\quad Q=\widehatP_1\widehatP_1^T, \quad \hat{C}=\widehatP_1^T\widehatP_0
\]
and 
\[
	U=(2P-I)(2Q-I).
\]

\begin{theorem}[Theorem~$5.2.2$ in~\cite{harmonyphd}]
\label{1-eignsp of U}
 	Let $P,Q$ be projections on $\mathbb{C}^m$.
	The $1$-eigenspace of $U$ is
	 \[ 
		\left(\col(P)\cap \col(Q)\right)\oplus \left(\ker(P)\cap \ker(Q)		\right)
	\] and it has dimension
 \[	
	m-\rk(P)-\rk(Q)+2\dim\left(\col(P)\cap \col(Q)\right).
\]
Moreover,
 \[
 	\col(P)\cap \col(Q)=\SPAN\{\mathbf{1}\}.
 \] 
\end{theorem}

\begin{theorem}[Lemma $2.3.6$ in~\cite{harmonyphd}]
\label{-1 eigen}
	The $(-1)$-eigenspace for $U$  is 
	\[
		\left(\col(P)\cap \ker(Q)\right)\oplus \left(\ker(P)\cap \col(Q)	\right)\] and its dimension is 
		\[
			\abs{C_0}+\abs{C_1}-2\rk(C).
			\]
\end{theorem}

\begin{theorem}[Lemma $2.3.7$ in \cite{harmonyphd}]
	Let $\mu\in(0,1)$ be an eigenvalue of $\hat{C}\hat{C}^T$. Choose $\theta	$ such that \[
	\cos\theta=2\mu-1.
	\]
	The map \[
	y\mapsto\left(\cos\theta+1\right)\widehatP_1y-\left(e^{i\theta}+1\right)	\widehatP_0\hat{C}^Ty
	\]
is an isomorphism from $\mu$-eigenspace of $\hat{C}\hat{C}^T$ to the $e^{i\theta}$-eigenspace of $U$, and the map
 \[
	y\mapsto\left(\cos\theta+1\right)\widehatP_1y-\left(e^{-i\theta}+1\right)\widehatP_0\hat{C}^Ty
\]
is an isomorphism from $\mu$-eigenspace of $\hat{C}\hat{C}^T$ to the $e^{-i\theta}$-eigenspace of $U$.
\end{theorem}

\begin{corollary}[Corollary $5.2.5$ in \cite{harmonyphd}]
\label{e^itheta-eigenspace of U}
Let $\mu\in(0,1)$ be an eigenvalue of $\hat{C}\hat{C}^T$. Choose $\theta$ such that $\cos\theta=2\mu-1.$ Let $E_\mu$ be the orthogonal projection onto the $\mu$-eigenspace of $\hat{C}\hat{C}^T$. Set 
\[
	W:= \widehatP_1E_\mu\widehatP_1^T. 
\]
Then the $e^{i\theta}$-eigenmatrix of $U$ is
\[
	\frac{1}{\sin^2(\theta)}\left( (\cos\theta+1)W-(e^{i\theta}+1)PW-(e^{-i		\theta}+1)WP+2PWP\right),
\] 
and the $e^{-i\theta}$-eigenmatrix of $U$ is
\[
	\frac{1}{\sin^2(\theta)}\left( (\cos\theta+1)W-(e^{-i\theta}+1)PW-(e^{i	\theta}+1)WP+2PWP\right).
\]
\end{corollary}

\section{Hamiltonians}

For every unitary matrix $U$, there exist Hermitian matrices $H$ such that
\[
	U=\exp(iH).
\] 
We call such $H$ a \textsl{Hamiltonian} of $U$. Since $U$ is unitary, it has spectral decomposition 
\[
	U=\sum_r e^{i\theta_r}E_r=\exp(iH),
\] 
and we can write
\[
	H=-i\sum_r \log(e^{i\theta_r})E_{\theta_r}=\sum_r \theta_r E_{\theta_r}.
\] 
For each eigenvalue $e^{i\theta_r}$ of $U$, we have that 
\[
	\log(e^{i\theta_r})=\log(e^{i\theta_r+2k_r\pi})
\] 
for non-zero integer $k_r$ and so, the choice of $H$ is not unique. That is, the Hamiltonian of $U$ is 
\[
	H=\sum_{\theta_r} (\theta_r+2k_r\pi) E_{\theta_r},
\] 
for any non-zero integer $k_r$. Note that $k_r$ are not necessarily equal for all the $\theta_r$.

Let $S$ be a real skew-symmetric matrix and $S$ can be viewed as the skew-adjacency matrix of a weighted oriented graph. When $H=iS$, we define the \textsl{$H$-digraph} to be the weighted oriented graph whose skew-adjacency matrix is $S$. This paper focuses on the case when the Hamiltonian can be written as $H=iS$ and studies the associated $H$-digraph.

For each eigenvalue $e^{i\theta_r}$ of $U$, if
$-\pi<\theta_r\leq \pi$ and $k_r=0$, the resulting unique Hamiltonian is called \textsl{principal Hamiltonian}. 
Let $H_0$ be the principle Hamiltonian. In general, if there is a real skew-symmetric $S_0$ such that $H_0=iS_0$, the choice 
\[
	H=H_0+\sum_r 2k_r\pi E_{\theta_r}
\] 
for non-constant $k_r$, cannot be written as $H=iS$ for a real skew-symmetric $S$.

Unless explicitly stated otherwise, we take the principal Hamiltonian
to be the Hamiltonian of $U$. Later in Corollary~\ref{A invertible}, we will show that there is a real skew-symmetric $S$ such that $H=iS$  if and only if the adjacency matrix of the bipartite graph $A(G)$ is invertible.  

\begin{theorem}
\label{A invertible theorem}
	Let $U$ be the transition matrix of the bipartite walk on a bipartite graph $G$. Let $H$ be the Hamiltonian of $U$ and let $E_{-1}$ be the projection onto the $(-1)$-eigenspace of $U$. Then there is a real skew-symmetric matrix $S$ such that
	\[
		H = iS +\pi E_{-1},
	\]
\end{theorem}

\proof 
Using the spectral decomposition  
\[
	U=\sum_r e^{i\theta_r}E_r=\exp(iH),
\] 
we can write
\[
	H=-i\sum_r \log(e^{i\theta_r})E_r=\sum_r \theta_r E_{\theta_r},
\] 
where $-\pi<\theta_r\leq \pi$. It follows that the $1$-eigenspace of $U$ corresponds to the $0$-eigenspace of $H$ and the $(-1)$-eigenspace of $U$ corresponds to the $\pi$-eigenspace of $H$ and $e^{i\theta_r}$-eigenspace gives $\theta_r$-eigenspace of $H$. 

Since $G$ is bipartite, the adjacency matrix of $G$ can be written as 
\[
	A(G)=
	\begin{pmatrix}
	\mathbf{0}&C\\
	C^T&\mathbf{0}
	\end{pmatrix}
\] 
for some $01$-matrix $C$. Let $\hat{C}$ be denoted the normalized version of $C$ and 
let $\mu\in(0,1)$ be an eigenvalue of $\hat{C}\hat{C}^T$. Choose $\theta$ such that
$\cos\theta=2\mu-1.$ Let $F_\mu$ be the orthogonal projection onto the $\mu$-eigenspace 
of $\hat{C}\hat{C}^T$. Set 
\[
	W:= \widehatP_1F_\mu\widehatP_1^T. 
\] 
By Corollary~\ref{e^itheta-eigenspace of U}, we have that 
\begin{align*}
H&=\sum_{\theta_r\neq \{1,-1\}}\theta_r\left(E_{\theta_r}-E_{-\theta_r}\right)+\pi\cdot E_{-1}\\
&=\sum_{\theta_r\neq \{1,-1\}}\theta_r\left(-\frac{2i}{\sin(\theta)}(PW-WP)\right)+\pi\cdot E_{-1}.
\end{align*}  
Since $\hat{C}\hat{C}^T$ is real and symmetric, we know that the orthogonal projection onto its $\mu$-eigenspace $F_\mu$ is real and symmetric. It follows that $W= \widehatP_1F_\mu\widehatP_1^T$ is real and symmetric. So the matrix $PW-WP$ is real. Set 
\[
	S=\sum_{\theta_r\neq \{1,-1\}}\theta_r\left(-\frac{2}{\sin(\theta_r)}(PW-WP)\right)
\]
and we know that $S$ is skew-symmetric.\qed

\begin{corollary}
\label{A invertible}
	Let $U$ be the transition matrix of the bipartite walk on a bipartite graph $G$. Let $S$ be a real skew-symmetric matrix and the Hamiltonian $H$ of $U$  can be written as $H=iS$ if and only if $A(G)$ is invertible.
\end{corollary}

\proof
By Theorem~\ref{-1 eigen}, we know that $E_{-1}$ is a real matrix. 
Using Theorem~\ref{A invertible theorem}, it is sufficient to prove that $E_{-1}=0$ if and only if $A(G)$ is invertible.

Now consider the $(-1)$-eigenvalue of $U$. From Theorem~\ref{-1 eigen} we know that 
\[
	\dim\left(E_{-1}\right)=\abs{C_0}+\abs{C_1}-2\rk(C).
\] 
This implies that $\dim\left(E_{-1}\right)=0$ if and only if 
\[
 \abs{C_0}+\abs{C_1}-2\rk(C)=0.
\]
Since $\rk(P_0)=\abs{C_0}$ and $\rk(P_1)=\abs{C_1}$ and $C=P_1^TP_0$, we get that  
\[
 \rk{C}\leq \min\{\abs{C_0},\abs{C_1}\}.
\] 
Thus, $\dim\left(E_{-1}\right)=0$ if and only if $\rk(P_0)=\rk(P_1)=\rk(C)$, which is equivalent to requiring that $C$ is invertible. Therefore we can conclude that there is a real skew-symmetric $S$ such that $H=iS$ if and only if $A(G)$ is invertible.\qed

Let $E_{\theta_r},E_{-\theta_r}$ be the corresponding eigenprojections  
of eigenvalue $e^{i\theta_r},e^{-i\theta_r}$ of $U$. Since $E_{\theta_r}$ are Hermitian, we have that 
\[
	E_{\theta_r}=\overline{E_{-\theta_r}}.
\] 
It follows that when $A(G)$ is invertible, the Hamiltonian 
\[
	H=\sum_{r} \theta_r \left( E_{\theta_r}-\overline{E_{\theta_r}}.\right)
\] 
has zero diagonal, which implies that the $H$-digraph has no loops.

We have proved that when $-1$ is an eigenvalue of $U$, there is no skew-symmetric matrix $S$ such 
that its Hamiltonian is in the form $H=iS$. So when $U$ has eigenvalue $-1$, 
we consider instead the Hamiltonian of $U^2$ and the $H$-digraph obtained from 
the Hamiltonian of $U^2$.

\section{Vertex-Face walks}
\label{vf walk}

Bipartite walks can be used to generalize many known walk models and one of them is the vertex-face walk. Here we show that vertex-face walk can be viewed as a special case of bipartite walk. As shown in~\cite{harmonyphd}, the Hamiltonian raised from vertex-face walk has many interesting properties, some of which will be presented using the bipartite walk language in this section and the next section.

An embedding of a graph $G$ in a surface $S$ is a continuous one-to-one map from $G$ to $S$. Given an embedding $G\rightarrow S$, the components of $S-G$ are called \textsl{regions}. If each region is homeomorphic to an open disk, then the embedding is called a \textsl{cellular embedding} and the regions are also called \textsl{faces} of the embedding. 

In \cite{harmonyphd}, Zhan introduces a new model of discrete quantum walk, the \textsl{vertex-face walk}. Let $\mathcal{M}$ be a circular embedding of graph $G$ on 
an orientable surface. Note that here the tail of the arc $(a,b)$ is vertex $a$. Let $M,N$ 
denote the arc-face incidence matrix and arc-tail incidence matrix respectively. The 
transition matrix of vertex-face walk on $\mathcal{M}$ is
\[
	U:= \left(2\widehat{M}\widehat{M}^T-I\right)\left(2\widehat{N}		\widehat{N}^T-I\right),
\]
where $\widehat{M},\widehat{N}$ is the matrices obtained from $M,N$ respectively by scaling each column to a unit vector.

The vertex-face incidence graph $X$ of the embedding $\mathcal{M}$ is a bipartite graph and two parts of $X$ are labelled by the vertices and the faces of $\mathcal{M}$. We can view the vertex-face walk on the circular embedding $\mathcal{M}$ as a bipartite walk by considering the bipartite walk over the vertex-face incidence graph of  $\mathcal{M}$. 

Now we show that the transition matrix of vertex-face walk on $\mathcal{M}$ is the same as the transition matrix of the bipartite walk on the vertex-face incidence graph of $\mathcal{M}$. Since $\mathcal{M}$ is a circular orientable embedding, the edges in the vertex-face incidence graph correspond to arcs of the embedding $\mathcal{M}$ of $G$. The arc-face incidence matrix $M$ of the embedding $\mathcal{M}$ is exactly the characteristic matrix of the edge-partition matrix  of the vertex-face incidence graph based on the face part. The arc-tail incidence matrix $N$ of the embedding $\mathcal{M}$ is exactly the characteristic matrix of the edge-partition matrix of the vertex-face incidence graph according to the vertex part. Hence, the bipartite walk on the incidence graph of the embedding $\mathcal{M}$ is exactly the same as the vertex-face walk on $\mathcal{M}$.
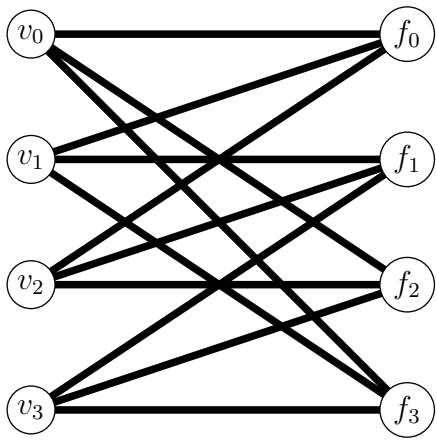
\begin{figure}[H]
    \centering
    \begin{subfigure}[t]{0.5\textwidth}
        \centering

\begin{tikzpicture}
\definecolor{cv0}{rgb}{0.0,0.0,0.0}
\definecolor{cfv0}{rgb}{1.0,1.0,1.0}
\definecolor{clv0}{rgb}{0.0,0.0,0.0}
\definecolor{cv1}{rgb}{0.0,0.0,0.0}
\definecolor{cfv1}{rgb}{1.0,1.0,1.0}
\definecolor{clv1}{rgb}{0.0,0.0,0.0}
\definecolor{cv2}{rgb}{0.0,0.0,0.0}
\definecolor{cfv2}{rgb}{1.0,1.0,1.0}
\definecolor{clv2}{rgb}{0.0,0.0,0.0}
\definecolor{cv3}{rgb}{0.0,0.0,0.0}
\definecolor{cfv3}{rgb}{1.0,1.0,1.0}
\definecolor{clv3}{rgb}{0.0,0.0,0.0}
\definecolor{cv0v1}{rgb}{0.0,0.0,0.0}
\definecolor{cv0v2}{rgb}{0.0,0.0,0.0}
\definecolor{cv0v3}{rgb}{0.0,0.0,0.0}
\definecolor{cv1v2}{rgb}{0.0,0.0,0.0}
\definecolor{cv1v3}{rgb}{0.0,0.0,0.0}
\definecolor{cv2v3}{rgb}{0.0,0.0,0.0}
\Vertex[style={minimum size=1.0cm,draw=cv0,fill=cfv0,text=clv0,shape=circle},LabelOut=false,L=\hbox{$1$},x=2.5cm,y=5.0cm]{v0}
\Vertex[style={minimum size=1.0cm,draw=cv1,fill=cfv1,text=clv1,shape=circle},LabelOut=false,L=\hbox{$2$},x=5.0cm,y=0cm]{v1}
\Vertex[style={minimum size=1.0cm,draw=cv2,fill=cfv2,text=clv2,shape=circle},LabelOut=false,L=\hbox{$3$},x=0cm,y=0.0cm]{v2}
\Vertex[style={minimum size=1.0cm,draw=cv3,fill=cfv3,text=clv3,shape=circle},LabelOut=false,L=\hbox{$0$},x=2.5cm,y=2cm]{v3}
\Edge[lw=0.1cm,style={color=cv0v1,},](v0)(v1)
\Edge[lw=0.1cm,style={color=cv0v2,},](v0)(v2)
\Edge[lw=0.1cm,style={color=cv0v3,},](v0)(v3)
\Edge[lw=0.1cm,style={color=cv1v2,},](v1)(v2)
\Edge[lw=0.1cm,style={color=cv1v3,},](v1)(v3)
\Edge[lw=0.1cm,style={color=cv2v3,},](v2)(v3)
\end{tikzpicture}
    
      \caption{The circular embedding of $K_4$}
    \end{subfigure}%
    ~
    
    The facial walks on $K_4$ embedding above: \begin{align*}
    f_0&=\{(0,1),(1,2),(2,0)\}\\
    f_1&=\{(1,3),(3,2),(2,1)\}\\
    f_2&=\{(0,2),(2,3),(3,0)\}\\
    f_3&=\{(0,3),(3,1),(1,0)\}
    \end{align*}
    
    ~\begin{subfigure}[t]{0.5\textwidth}
        \centering
   
\begin{tikzpicture}[rotate=270]
\definecolor{cv0}{rgb}{0.0,0.0,0.0}
\definecolor{cfv0}{rgb}{1.0,1.0,1.0}
\definecolor{clv0}{rgb}{0.0,0.0,0.0}
\definecolor{cv1}{rgb}{0.0,0.0,0.0}
\definecolor{cfv1}{rgb}{1.0,1.0,1.0}
\definecolor{clv1}{rgb}{0.0,0.0,0.0}
\definecolor{cv2}{rgb}{0.0,0.0,0.0}
\definecolor{cfv2}{rgb}{1.0,1.0,1.0}
\definecolor{clv2}{rgb}{0.0,0.0,0.0}
\definecolor{cv3}{rgb}{0.0,0.0,0.0}
\definecolor{cfv3}{rgb}{1.0,1.0,1.0}
\definecolor{clv3}{rgb}{0.0,0.0,0.0}
\definecolor{cv4}{rgb}{0.0,0.0,0.0}
\definecolor{cfv4}{rgb}{1.0,1.0,1.0}
\definecolor{clv4}{rgb}{0.0,0.0,0.0}
\definecolor{cv5}{rgb}{0.0,0.0,0.0}
\definecolor{cfv5}{rgb}{1.0,1.0,1.0}
\definecolor{clv5}{rgb}{0.0,0.0,0.0}
\definecolor{cv6}{rgb}{0.0,0.0,0.0}
\definecolor{cfv6}{rgb}{1.0,1.0,1.0}
\definecolor{clv6}{rgb}{0.0,0.0,0.0}
\definecolor{cv7}{rgb}{0.0,0.0,0.0}
\definecolor{cfv7}{rgb}{1.0,1.0,1.0}
\definecolor{clv7}{rgb}{0.0,0.0,0.0}
\definecolor{cv0v4}{rgb}{0.0,0.0,0.0}
\definecolor{cv0v5}{rgb}{0.0,0.0,0.0}
\definecolor{cv0v6}{rgb}{0.0,0.0,0.0}
\definecolor{cv0v7}{rgb}{0.0,0.0,0.0}
\definecolor{cv1v4}{rgb}{0.0,0.0,0.0}
\definecolor{cv1v5}{rgb}{0.0,0.0,0.0}
\definecolor{cv1v6}{rgb}{0.0,0.0,0.0}
\definecolor{cv1v7}{rgb}{0.0,0.0,0.0}
\definecolor{cv2v4}{rgb}{0.0,0.0,0.0}
\definecolor{cv2v5}{rgb}{0.0,0.0,0.0}
\definecolor{cv2v6}{rgb}{0.0,0.0,0.0}
\definecolor{cv2v7}{rgb}{0.0,0.0,0.0}
\definecolor{cv3v4}{rgb}{0.0,0.0,0.0}
\definecolor{cv3v5}{rgb}{0.0,0.0,0.0}
\definecolor{cv3v6}{rgb}{0.0,0.0,0.0}
\definecolor{cv3v7}{rgb}{0.0,0.0,0.0}
\Vertex[style={minimum size=1.0cm,draw=cv0,fill=cfv0,text=clv0,shape=circle},LabelOut=false,L=\hbox{$f_0$},x=0.0cm,y=5.0cm]{v0}
\Vertex[style={minimum size=1.0cm,draw=cv1,fill=cfv1,text=clv1,shape=circle},LabelOut=false,L=\hbox{$f_1$},x=1.6667cm,y=5.0cm]{v1}
\Vertex[style={minimum size=1.0cm,draw=cv2,fill=cfv2,text=clv2,shape=circle},LabelOut=false,L=\hbox{$f_2$},x=3.3333cm,y=5.0cm]{v2}
\Vertex[style={minimum size=1.0cm,draw=cv3,fill=cfv3,text=clv3,shape=circle},LabelOut=false,L=\hbox{$f_3$},x=5.0cm,y=5.0cm]{v3}
\Vertex[style={minimum size=1.0cm,draw=cv4,fill=cfv4,text=clv4,shape=circle},LabelOut=false,L=\hbox{$v_0$},x=0.0cm,y=0.0cm]{v4}
\Vertex[style={minimum size=1.0cm,draw=cv5,fill=cfv5,text=clv5,shape=circle},LabelOut=false,L=\hbox{$v_1$},x=1.6667cm,y=0.0cm]{v5}
\Vertex[style={minimum size=1.0cm,draw=cv6,fill=cfv6,text=clv6,shape=circle},LabelOut=false,L=\hbox{$v_2$},x=3.3333cm,y=0.0cm]{v6}
\Vertex[style={minimum size=1.0cm,draw=cv7,fill=cfv7,text=clv7,shape=circle},LabelOut=false,L=\hbox{$v_3$},x=5.0cm,y=0.0cm]{v7}
\Edge[lw=0.1cm,style={color=cv0v4,},](v0)(v4)
\Edge[lw=0.1cm,style={color=cv0v5,},](v0)(v5)
\Edge[lw=0.1cm,style={color=cv0v6,},](v0)(v6)
\Edge[lw=0.1cm,style={color=cv1v5,},](v1)(v5)
\Edge[lw=0.1cm,style={color=cv1v6,},](v1)(v6)
\Edge[lw=0.1cm,style={color=cv1v7,},](v1)(v7)
\Edge[lw=0.1cm,style={color=cv2v4,},](v2)(v4)
\Edge[lw=0.1cm,style={color=cv2v6,},](v2)(v6)
\Edge[lw=0.1cm,style={color=cv2v7,},](v2)(v7)
\Edge[lw=0.1cm,style={color=cv3v4,},](v3)(v4)
\Edge[lw=0.1cm,style={color=cv3v5,},](v3)(v5)
\Edge[lw=0.1cm,style={color=cv3v7,},](v3)(v7)
\end{tikzpicture}
\caption{ The vertex-face incidence graph of the planar embedding of $K_4$}
    \end{subfigure}
    \caption{The circular embedding of $K_4$ and its corresponding vertex-face incidence graph}
\end{figure}

In \cite{harmonyphd}, Zhan focuses on the circular orientable embedding of graph $G$ such that both $G$ and its dual graph are regular. The embedding $\mathcal{M}$ has type $(k,l)$ if each vertex has degree $l$ and each faces uses $k$ vertices. Note that a vertex-face walk over a $(k,l)$-type embedding $\mathcal{M}$ corresponds to a bipartite walk on a $(k,l)$-regular bipartite graph that is the vertex-face incidence graph of $\mathcal{M}$.

\begin{theorem}[Theorem~$8.5.4$ in~\cite{DQW}]
Let $G$ be a semi-regular bipartite graph with degree $(k,l)$ and $P_0,P_1$ denote its two parts. Let $\pi_0,\pi_1$ denote the partitions of edges of $G$ according to $P_0,P_1$ respectively. Let $U$ be the bipartite walk transition matrix for $G$. Then 
\[
	U^2=\exp\left(\gamma(U-U^T)\right)
\] for some real number $\gamma$ if and only if $G$ has four or five distinct eigenvalues. Moreover, 
\[ S=
	\frac{kl}{4}(U^T-U)
\] is the skew-adjacency matrix of some oriented graph on the edges of $G$.

 Let $c_{0,k}$ denote the cell of partition $\pi_0$ containing edge $e_k$ and similarly, $c_{1,k}$ denote the cell of partition $\pi_1$ containing edge $e_k$. Then we have 
 \[
	S_{i,j}=\begin{cases}
 	1,\quad\text{if }\abs{c_{0,i}\cap c_{1,j}}=1\text{ and }\abs{c_{1,i}\cap c_{0,j}}=0, \\[2.5mm]
 	-1,\quad\text{if }\abs{c_{0,i}\cap c_{1,j}}=0\text{ and }\abs{c_{1,i}\cap c_{0,j}}=1, \\[2.5mm]
 	0,\quad\text{otherwise.}
 	\end{cases}\qed
 	\]
\end{theorem}

A partial geometric design with parameters $(d,k,t,c)$ is a point-$d$-regular and block-$k$-regular design, where for each point-block pair $(p,B)$, the number of incident point-block pairs\[
\abs{ \{(p',B'): p'\neq p, B'\neq B, p'\in B,p\in B'\}} \]equals $c$ or $t$, depending on whether $p$ is in $B$ or not.
In~\cite{DQW} Theorem~$8.5.5$, Godsil and Zhan have showed that when $G$ is an incidence graph of a partial geometric design, then we have that
\[
	U^2=\exp\left(\gamma(U-U^T)\right)
\] for some real number $\gamma$.

\section{Vertex-Face walks on complete graphs}
\label{vxf on CompleteG}
In \cite{biggembedding}, Biggs states that $K_n$ has a regular embedding if and only if $n$ is a prime power and every regular embedding of $K_n$ must arise from the rotation system stated in~\cite{harmonyphd}.

\begin{lemma}[Theorem $5.6.2$ in~\cite{harmonyphd}] Let $n=p^k$ for some prime $p$. Let $g$ be a primitive generator of the finite field $\mathbb{F}$ of order $n$. For each element $u$ in $\mathbb{F}$, define the cyclic permutation 
\[
	\pi_u=\{v+g^0,v+g^1,\cdots,v+g^{n-2}\}.
	\]
The rotation system $\{\pi_u:u\in V(K_m)\}$ gives a circular embedding of $K_n$.
\end{lemma} 

In the case of $H$-digraphs arised from the vertex-face walk on $K_n$, we know that the skew-adjacency matrix of $H$-digraph $A\big(\overrightarrow{H}\big)$ is indexed by arcs of $K_n$. 
 Let $f_{ab}$ denote the unique face that contains arc $(a,b)$.
From the proof of Theorem~$8.5.4$ in \cite{DQW}, we have that 
\[
 A\big(\overrightarrow{H}\big)_{(a,b),(c,d)}=\begin{cases}
 1,\quad\text{if }c\in f_{ab}\text{ and }a\not\in f_{cd}, \\[2.5mm]
 -1,\quad\text{if }a\in f_{cd}\text{ and }c\not\in f_{ab}, \\[2.5mm]
 0,\quad\text{otherwise.}
 \end{cases}
 \]
 Note that in a self-dual circular embedding of $K_n$, each face consists of $n-1$ distinct vertices, which implies that each face misses a unique vertex of $K_n$.
 
We use $LD\left(K_n\right)$ to denote the line digraph of $K_n$.
 
\begin{theorem}
\label{H-digraph of vx-fc of K_n}
The $H$-digraphs $Z_n$ obtained from the vertex-face walks of a self-dual embedding of $K_n$ 
is the line digraphs of $K_n$.
\end{theorem}

\proof 
We construct an isomorphism from $Z_n$ to $LD(K_n)$. Define a map $f:V(Z_n)\rightarrow V\left(LD(K_n)\right)$ as
\[
	(a,b)\mapsto (u,a),
\] 
where $u$ is the unique vertex missed by $f_{ab}$. First we show that $f$ is a homomorphism. Say
\[
	f(a,b)=(u,a),\quad f(c,d)=(v,c),
\] 
which implies that $u$ is the unique vertex missed by $f_{ab}$ and $v$ is the unique vertex missed by $f_{cd}$. We know that there is an arc from $(a,b)$ to $(c,d)$ in $Z_n$ if and only if 
\[
	c\in f_{ab}\text{ and }a\not\in f_{cd}.
\] 
Since each face miss a unique vertex in the circular embedding of $K_n$, we must have that 
\[
	a=v,
\] which means that there is an arc from $f(a,b)$ to $f(c,d)$ in $LD(K_n)$. Thus, the map $f$ is indeed a homomorphism.

Now we prove that $f$ is a bijection and since $LD(K_n)$ is finite, it suffices to prove that $f$ is an injection. Assume towards contradictions that two distinct arcs $(a,b)$ and $(a',b')$ get mapped to $(x,y)$ by the map $f$. Then by how we define the map $f$, we know that
\[
	a=a'=y.
\] 
The vertex $x$ is missed by $f_{ab}$ and $f_{a'b'}=f_{ab'}$. Since the faces here arised from facial walks on the circular embedding of $K_n$, we must have that 
\[
	(a,b)=(a',b').
\] 
This means that $f$ has to be an injection and hence, a bijection. Therefore, we can conclude that the map $f$ gives an isomorphism from $Z_n$ to $LD(K_n)$. \qed

\begin{theorem}[Theorem~$5.6.3$ in \cite{harmonyphd}]
Let $n$ be a prime power. Let $U$ be the transition matrix of the vertex-face walk for a regular embedding of $K_n$. Then there is a $\gamma\in\mathbb{R}$ such that
\[
	U=\exp\left(\gamma(U^T-U)\right).
\] 
Further $U^T-U$ is a scalar multiple of the skew-adjacency matrix of an oriented graph, which 
\begin{enumerate}[label=(\roman*)]
\item has $n(n-1)$ vertices,
\item is $(n-2)$-regular, and
\item has exactly three eigenvalues: $0$ and $\pm i\sqrt{n(n-2)}$
\end{enumerate}
\end{theorem}

We rephrase Theorem~\ref{H-digraph of vx-fc of K_n} in terms of bipartite walk and we get the following theorem.
 
\begin{theorem}
Let $G_n$ be a $(n-1)$-regular bipartite graph with each part of size $n$. Then the $H$-digraph obtained from the bipartite walk on $G_n$ is the line digraph of $K_n$.
\end{theorem}

\proof 
Since there is every cell of $\pi_1$ miss a unique vertex in $C_0$ and every cell of $\pi_0$ misses a unique vertex in $C_1$, the proof of Theorem~\ref{H-digraph of vx-fc of K_n} applies here.\qed

\section{Paths and even cycles}
The vertex-face incidence graph of a cellular embedding of a graph must have degree at least three for each vertex. So neither a path nor a cycle can be a bipartite graph raised from the vertex-face incidence relation of an circular embedding. In this section, we discuss the bipartite walk defined on paths and even cycles.
\label{bipartite walk on paths}

\begin{figure}[H]
\begin{center}
\begin{tikzpicture}[scale=0.55]
\definecolor{cv0}{rgb}{0.0,0.0,0.0}
\definecolor{cfv0}{rgb}{1.0,1.0,1.0}
\definecolor{clv0}{rgb}{0.0,0.0,0.0}
\definecolor{cv1}{rgb}{0.0,0.0,0.0}
\definecolor{cfv1}{rgb}{1.0,1.0,1.0}
\definecolor{clv1}{rgb}{0.0,0.0,0.0}
\definecolor{cv2}{rgb}{0.0,0.0,0.0}
\definecolor{cfv2}{rgb}{1.0,1.0,1.0}
\definecolor{clv2}{rgb}{0.0,0.0,0.0}
\definecolor{cv3}{rgb}{0.0,0.0,0.0}
\definecolor{cfv3}{rgb}{1.0,1.0,1.0}
\definecolor{clv3}{rgb}{0.0,0.0,0.0}
\definecolor{cv4}{rgb}{0.0,0.0,0.0}
\definecolor{cfv4}{rgb}{1.0,1.0,1.0}
\definecolor{clv4}{rgb}{0.0,0.0,0.0}
\definecolor{cv5}{rgb}{0.0,0.0,0.0}
\definecolor{cfv5}{rgb}{1.0,1.0,1.0}
\definecolor{clv5}{rgb}{0.0,0.0,0.0}
\definecolor{cv6}{rgb}{0.0,0.0,0.0}
\definecolor{cfv6}{rgb}{1.0,1.0,1.0}
\definecolor{clv6}{rgb}{0.0,0.0,0.0}
\definecolor{cv7}{rgb}{0.0,0.0,0.0}
\definecolor{cfv7}{rgb}{1.0,1.0,1.0}
\definecolor{clv7}{rgb}{0.0,0.0,0.0}
\definecolor{cv0v1}{rgb}{0.0,0.0,0.0}
\definecolor{cv1v2}{rgb}{0.0,0.0,0.0}
\definecolor{cv2v3}{rgb}{0.0,0.0,0.0}
\definecolor{cv3v4}{rgb}{0.0,0.0,0.0}
\definecolor{cv4v5}{rgb}{0.0,0.0,0.0}
\definecolor{cv5v6}{rgb}{0.0,0.0,0.0}
\definecolor{cv6v7}{rgb}{0.0,0.0,0.0}
\Vertex[style={minimum size=1.0cm,draw=cv0,fill=cfv0,text=clv0,shape=circle},LabelOut=false,L=\hbox{$0$},x=4.5cm,y=5cm]{v0}
\Vertex[style={minimum size=1.0cm,draw=cv1,fill=cfv1,text=clv1,shape=circle},LabelOut=false,L=\hbox{$1$},x=0cm,y=5cm]{v1}
\Vertex[style={minimum size=1.0cm,draw=cv2,fill=cfv2,text=clv2,shape=circle},LabelOut=false,L=\hbox{$2$},x=4.5cm,y=3cm]{v2}
\Vertex[style={minimum size=1.0cm,draw=cv3,fill=cfv3,text=clv3,shape=circle},LabelOut=false,L=\hbox{$3$},x=0cm,y=3cm]{v3}
\Vertex[style={minimum size=1.0cm,draw=cv4,fill=cfv4,text=clv4,shape=circle},LabelOut=false,L=\hbox{$4$},x=4.5cm,y=1cm]{v4}
\Vertex[style={minimum size=1.0cm,draw=cv5,fill=cfv5,text=clv5,shape=circle},LabelOut=false,L=\hbox{$5$},x=0cm,y=1cm]{v5}
\Vertex[style={minimum size=1.0cm,draw=cv6,fill=cfv6,text=clv6,shape=circle},LabelOut=false,L=\hbox{$6$},x=4.5cm,y=-1cm]{v6}
\Vertex[style={minimum size=1.0cm,draw=cv7,fill=cfv7,text=clv7,shape=circle},LabelOut=false,L=\hbox{$7$},x=0cm,y=-1cm]{v7}
\Edge[lw=0.01cm,style={color=cv0v1,},label={$e_0$}](v0)(v1)
\Edge[lw=0.01cm,style={color=cv1v2,},label={$e_1$}](v1)(v2)
\Edge[lw=0.01cm,style={color=cv2v3,},label={$e_2$}](v2)(v3)
\Edge[lw=0.01cm,style={color=cv3v4,},label={$e_3$}](v3)(v4)
\Edge[lw=0.01cm,style={color=cv4v5,},label={$e_4$}](v4)(v5)
\Edge[lw=0.01cm,style={color=cv5v6,},label={$e_5$}](v5)(v6)
\Edge[lw=0.01cm,style={color=cv6v7,},label={$e_6$}](v6)(v7)
\end{tikzpicture}
\end{center}
\caption{$P_8$}
\end{figure}

We label the vertices of $P_n$ as $v_0,v_1\cdots,v_{n-1}$ accordingly from the leftmost vertices to the rightmost vertices of $P_n$. Note that $v_0,v_{n-1}$ are the only two vertices of degree $1$ with all the others of degree $2$. Partition $\pi_0$ is the partition of edges such that edges with the same end at a vertex in $\{v_1,v_3,\cdots,v_{n-1}\}$ are in the same cell of $\pi_0$. Partition $\pi_1$ is the partition of edges such that edges with the same end at a vertex in $ \{v_0,v_2,\cdots,\allowbreak v_{n-2}\}$  are in the same cell of $\pi_1$. Edge $e_i$ is  the edge between $v_i,v_{i+1}$ for all integer $0\leq i\leq n-2$.

Recall that $P,Q$ are the projections onto the vectors that is constant on the cells of $\pi_0,\pi_1$ respectively. Let $c_i$ denote the characteristic vector of the edges adjacent to vertex $i$. The column space of $Q$ is 
 \[
	\col(Q)=\SPAN\{c_0,c_2,\cdots,c_{n-2}\},
	\]
The matrix $2Q-I$ is a reflection about the column space of $Q$, which is the span of cells of $\pi_1$. If two edges belong to the same cell, then they are the ``cellmate" of each other.

Note that every vertex of a path has degree $\leq 2$, which means that each edge has at most one cellmate in the partitions. For each $0\leq i\leq n-2$, let $e_j$ be the cellmate of $e_i$ in $\pi_1$. Using that each cell in $\pi_0,\pi_1$ has size $\leq 2$, we have that 
\[
	(2Q-I)e_i=e_j.
\] Similarly, if $e_i,e_j$ are cellmates in $\pi_0$, then we have that 
\[
	(2P-I)e_i=e_j.
\] Here both reflections $2P-I$ and $2Q-I$ is permutation matrices.
Thus, the transition matrix $U=(2P-I)(2Q-I)$ of bipartite walk on $P_n$
is a permutation matrix such that for each integer $0\leq i\leq n-2$,  \begin{equation}
\label{path U permutation}
Ue_i=
\begin{cases}
e_{i+2}, \quad\text{if }i \text{ is odd and }i\neq n-3;\\
e_{i-2}, \quad\text{if }i \text{ is even and }i\neq 0; \\
e_{1}, \quad\text{if }i =0;\\
e_{n-2},\quad\text{if }i =n-3.
\end{cases}
\end{equation}

\begin{theorem}
\label{path U permutation}
The transition matrix of the bipartite walk on $P_n$ corresponds to a $(n-1)$-cycle permutation whose cycle form is 
\[
	\left(e_0,e_1,e_3,\cdots,e_{n-3},e_{n-2},e_{n-4},\cdots,e_2\right).
\]
\end{theorem}

\proof It follows from the discussion above.\qed

 For example, the transition matrix of the bipartite walk on $P_8$ is 
\[
	U=
\begin{pmatrix}
0 & 0 & 1 & 0 & 0 & 0& 0 \\
1 & 0 & 0 & 0 & 0 & 0 & 0 \\
0 & 0 & 0 & 0 & 1& 0 & 0 \\
0 & 1& 0 & 0 & 0 & 0 & 0 \\
0 & 0 & 0 & 0 & 0 & 0 & 1 \\
0 & 0 & 0 & 1&0 & 0 & 0 \\
0 & 0 & 0 & 0 & 0& 1 & 0
\end{pmatrix}.
\] This correspond to the permutation $(0135642)$ in $S_7$  and we have that \[
	U^7=I.
	\] Since $U(P_8)$ is a permutation matrix of order $7$, it is easy to see that every edge of $P_8$ can be mapped to any other edges within $7$ steps in the bipartite walk. This is an interesting phenomenon called \textsl{universal perfect state transfer}. Note that if $U$ is the transition matrix of bipartite walk on $P_n$, then 
\[
	U^{n-1}=I,
\]which implies that for every $n$, the bipartite walk on $P_n$ has the universal perfect state transfer. We will discuss this property further in the next section. 

Cyclic permutation matrix $U$ is of order $n-1$, then it has eigenvalue
\[
	\lambda_k=\left(e^{\frac{2\pi i}{n-1}}\right)^k
\] with eigenvector 
\begin{equation}
\label{eigenvector of U(Pn)}
f_k=\begin{pmatrix}
1&
\lambda_k^{-1}&
\lambda_k&
\lambda_k^{-2}&
\lambda_k^{2}&
\cdots&
\lambda_k^{-(n-2)/2}&
\lambda_k^{(n-2)/2}
\end{pmatrix}^T,
\end{equation}
for $k=0,\cdots,n-2$. 
The $\lambda_k$-eigenspace of $U$ is 
\[
	E_{\lambda_k}=\frac{1}{n-1}ff^*.
\] 
Note that $E_1=\frac{1}{n-1}J.$

From the eigenvectors of $U$~(\ref{eigenvector of U(Pn)}), we know that if $s,t$ are integers in $\{1,\cdots,n-2\}$, we have that 
\begin{equation}
\label{entry of E_r of path}
\left(E_{\lambda_r}\right)_{s,t}=
	\begin{cases}\frac{1}{n-1}(\lambda_r)^{-\frac{s+1}{2}}(\lambda_r)^{\frac{t+1}{2}}\quad\text{if both }s, t \text{ are odd;}\\[2mm]

\frac{1}{n-1}(\lambda_r)^{\frac{s}{2}}(\lambda_r)^{\frac{t+1}{2}}\quad
\text{if }s\text{ is even and } t \text{ is odd;}\\[2mm]

\frac{1}{n-1}(\lambda_r)^{-\frac{s+1}{2}}(\lambda_r)^{-\frac{t}{2}}\quad\text{if }s\text{ is odd and } t \text{ is even;}\\[2mm]

\frac{1}{n-1}(\lambda_r)^{\frac{s}{2}}(\lambda_r)^{-\frac{t}{2}}\quad\text{if both }s, t \text{ are even.}
\end{cases}
\end{equation}

\begin{theorem}
\label{H-digraph of Path}
For an even $n\geq 4$, the $H$-digraph obtained from the bipartite walk on $P_n$ is an oriented $K_{n-1}$.
\end{theorem}

\proof
As the discussion above, the transition matrix of bipartite walk on $P_n$ has spectral decomposition
\[
	U=\sum_{k=0}^{n-2} \lambda_k E_{\lambda_k},
	\] where 
\[	
	\lambda_k=\left(e^{\frac{2\pi i}{n-1}}\right)^k.
\]
When $n$ is even, the Hamiltonian of $U$ is 
\[
	H=\sum_{k=0}^{(n-2)/2} \frac{2k\pi }{n-1}\left(E_{\lambda_k}-
\overline{E_{\lambda_k}}\right).
\]

To prove that the $H$-digraph is an oriented complete graph, we  show that the Hamiltonian $H$ has non-zero off-diagonal entries. As shown above that the eigenvector of $U$ with eigenvalue $\lambda_k$ is of the form~\ref{eigenvector of U(Pn)}, each row of $E_{\lambda_k}$ is a permutation of its first row, which implies that each row of $H$ is a permutation of its first row. So in order to prove that all the off-diagonal entries of $H$ are non-zero, it is sufficient to prove that 
\[
	H_{0,t}\neq 0
\] for all $t\neq 0$.

Based on the formula of the $(s,t)$-th entry of $E_{\lambda_r}$ shown in~\ref{entry of E_r of path}
we have that
for $r\in\{0,1,2,\cdots,n-2\} $ and, $s,t\in\{0,1,\cdots,n-2\}$, we have that 
\[
\left(E_{\lambda_r}-\overline{E_{\lambda_r}}\right)_{s,t}=\begin{cases}
\frac{2}{n-1}\sin\left(\frac{2\pi r}{n-1}\cdot \frac{t-s}{2}\right)i,\quad\text{if both }s, t \text{ are odd;}\\[2mm]

\frac{2}{n-1}\sin\left(\frac{2\pi r}{n-1}\cdot\frac{s+t+1}{2} \right)i \quad\text{if }s\text{ is even and } t \text{ is odd;}\\[2mm]

\frac{2}{n-1}\sin\left(\frac{2\pi r}{n-1}\cdot \frac{-t-s-1}{2}\right)i,\quad\text{if }s\text{ is odd and } t \text{ is even;}\\[2mm]
\frac{2}{n-1}\sin\left(\frac{2\pi r}{n-1}\cdot \frac{s-t}{2}\right)i,\quad\text{if both }s, t \text{ are even.}\\[2mm]
0\quad\text{if }s=t.
\end{cases}
\] 
Then entries of the first row of $H$ are 
\[
\left( H\right)_{0,t}=\sum_{k=0}^{(n-2)/2} \frac{2k\pi }{n-1}\left(E_{\lambda_k}\right)_{0,t}= \begin{cases}
\sum_{k=0}^{(n-2)/2} \frac{4k\pi }{(n-1)^2}\sin\left(\frac{2k\pi}{n-1}\cdot\frac{t+1}{2}\right) \quad\text{if }t \text{ is odd;}\\[3mm]
\sum_{k=0}^{(n-2)/2} \frac{4k\pi }{(n-1)^2}\sin\left(\frac{2k\pi}{n-1}\cdot\frac{-t}{2}\right) \quad\text{if }t \text{ is even;}\\[3mm]
0  \quad\text{if }t=0
\end{cases}
\]
When $n=2a+2$ for some integer $a\geq 1$, then for each positive odd integer $b$, we have that
\begin{equation}
\label{sum sine}
\sum_{k=0}^{(n-2)/2} \frac{2k\pi }{n-1}\sin\left(\frac{2k\pi}{n-1}\cdot b\right)
=\frac{\pi\csc\left(\frac{b\pi}{2a+1}\right)\left(2(a+1)+\sin\left(\frac{2b\pi(a+1)}{2a+1}\right)\csc\left(\frac{b\cdot\pi}{2a+1}\right)\right)}{4a+2}
\end{equation} and
for each positive even integer $b$, we have that
\begin{equation}
\label{even sum sine}
\sum_{k=0}^{(n-2)/2} \frac{2k\pi }{n-1}\sin\left(\frac{2k\pi}{n-1}\cdot b\right)
=\frac{\pi\csc\left(\frac{b\pi}{2a+1}\right)\left(-2(a+1)+\sin\left(\frac{2b\pi(a+1)}{2a+1}\right)\csc\left(\frac{b\cdot\pi}{2a+1}\right)\right)}{4a+2}.
\end{equation}
Since the sine function is an odd function, we only need to show that $H_{0,t}\neq 0$ for all odd $1\leq t\leq \frac{n}{2}$.
Since $\csc(x)\neq 0$ over all its domain and when $1\leq b\leq a+1 $,
\[ 
	\sin\left(\frac{2b\pi(a+1)}{2a+1}\right)\csc\left(\frac{b\cdot\pi}{2a+1}\right)\pm 2(a+1)\neq 0.
\]
The sum shown in~\ref{even sum sine} and \ref{sum sine} are non-zero for all $1\leq b\leq a+1$. Thus, we have that 
\[
	\left( H\right)_{0,t}\neq 0
\] 
for all $t\neq 0$.
Therefore, we can conclude that the $H$-digraph is an oriented $K_{n-1}$.\qed

Note that when $n$ is odd, the adjacency matrix of $P_n$ is not invertible and so we consider the Hamiltonian of $U^2$. When $n=3$, the Hamiltonian of $U^2$ is zero matrix. When $n\equiv 1\Mod 4$, the square of its transition matrix $U^2$ still has $-1$ as an eigenvalue, which implies that there is no real skew-symmetric $S$ such that Hamiltonian of $U^2$ is of the form $iS$. So here, we omit the case when $n\equiv 1\Mod 4$.

\begin{corollary}
\label{odd path H-digraph}
When $n\equiv 3\Mod 4$, let
\[
 U^2=\exp(iH),
 \] 
 then $H$ is the weighted skew adjacency matrix of two copies of oriented $K_{\frac{n-1}{2}}$. 
\end{corollary}

\proof By Theorem~\ref{path U permutation}, we know that $U^2$ corresponds to two $\left(\frac{n-1}{2}\right)$-cycles. Each $\left(\frac{n-1}{2}\right)$-cycle is equivalent to the permutation associated with the transition matrix of $P_{\frac{n+1}{2}}$. The result follows from Theorem~\ref{H-digraph of Path}.\qed

Even cycles are another class of bipartite graphs that cannot be raised from the vertex-face incidence relation of a circular embedding.

For an even integer $n$, consider a path $P_n$ with the same labelling as before and add an edge $e_{n-1}$ between $v_0,v_{n-1}$, which gives us a even cycle $C_n$.
Partition $\pi_0$ are the partition of edges based on vertices $\{v_1,v_3,\cdots,v_{n-1}\}$ and partition $\pi_1$ are the partition of edges based on vertices $\{v_0,v_2,\cdots,v_{n-2}\}$ . 
\begin{figure}[H]
\begin{center}
\begin{tikzpicture}[scale=0.7]
\definecolor{cv0}{rgb}{0.0,0.0,0.0}
\definecolor{cfv0}{rgb}{1.0,1.0,1.0}
\definecolor{clv0}{rgb}{0.0,0.0,0.0}
\definecolor{cv1}{rgb}{0.0,0.0,0.0}
\definecolor{cfv1}{rgb}{1.0,1.0,1.0}
\definecolor{clv1}{rgb}{0.0,0.0,0.0}
\definecolor{cv2}{rgb}{0.0,0.0,0.0}
\definecolor{cfv2}{rgb}{1.0,1.0,1.0}
\definecolor{clv2}{rgb}{0.0,0.0,0.0}
\definecolor{cv3}{rgb}{0.0,0.0,0.0}
\definecolor{cfv3}{rgb}{1.0,1.0,1.0}
\definecolor{clv3}{rgb}{0.0,0.0,0.0}
\definecolor{cv4}{rgb}{0.0,0.0,0.0}
\definecolor{cfv4}{rgb}{1.0,1.0,1.0}
\definecolor{clv4}{rgb}{0.0,0.0,0.0}
\definecolor{cv5}{rgb}{0.0,0.0,0.0}
\definecolor{cfv5}{rgb}{1.0,1.0,1.0}
\definecolor{clv5}{rgb}{0.0,0.0,0.0}
\definecolor{cv6}{rgb}{0.0,0.0,0.0}
\definecolor{cfv6}{rgb}{1.0,1.0,1.0}
\definecolor{clv6}{rgb}{0.0,0.0,0.0}
\definecolor{cv7}{rgb}{0.0,0.0,0.0}
\definecolor{cfv7}{rgb}{1.0,1.0,1.0}
\definecolor{clv7}{rgb}{0.0,0.0,0.0}
\definecolor{cv0v1}{rgb}{0.0,0.0,0.0}
\definecolor{cv1v2}{rgb}{0.0,0.0,0.0}
\definecolor{cv2v3}{rgb}{0.0,0.0,0.0}
\definecolor{cv3v4}{rgb}{0.0,0.0,0.0}
\definecolor{cv4v5}{rgb}{0.0,0.0,0.0}
\definecolor{cv5v6}{rgb}{0.0,0.0,0.0}
\definecolor{cv6v7}{rgb}{0.0,0.0,0.0}
\Vertex[style={minimum size=1.0cm,draw=cv0,fill=cfv0,text=clv0,shape=circle},LabelOut=false,L=\hbox{$0$},x=4.5cm,y=6cm]{v0}
\Vertex[style={minimum size=1.0cm,draw=cv1,fill=cfv1,text=clv1,shape=circle},LabelOut=false,L=\hbox{$1$},x=0cm,y=6cm]{v1}
\Vertex[style={minimum size=1.0cm,draw=cv2,fill=cfv2,text=clv2,shape=circle},LabelOut=false,L=\hbox{$2$},x=4.5cm,y=2.5cm]{v2}
\Vertex[style={minimum size=1.0cm,draw=cv3,fill=cfv3,text=clv3,shape=circle},LabelOut=false,L=\hbox{$3$},x=0cm,y=2.5cm]{v3}
\Vertex[style={minimum size=1.0cm,draw=cv4,fill=cfv4,text=clv4,shape=circle},LabelOut=false,L=\hbox{$4$},x=4.5cm,y=0.5cm]{v4}
\Vertex[style={minimum size=1.0cm,draw=cv5,fill=cfv5,text=clv5,shape=circle},LabelOut=false,L=\hbox{$5$},x=0cm,y=0.5cm]{v5}
\Edge[lw=0.01cm,style={color=cv0v1,},label={$e_0$}](v0)(v1)
\Edge[lw=0.01cm,style={color=cv5v6,},label={$e_1$}](v1)(v2)
\Edge[lw=0.01cm,style={color=cv1v2,},label={$e_2$}](v2)(v3)
\Edge[lw=0.01cm,style={color=cv2v3,},label={$e_3$}](v3)(v4)
\Edge[lw=0.01cm,style={color=cv3v4,},label={$e_4$}](v5)(v4)
\Edge[lw=0.01cm,style={color=cv4v5,},label={$e_5$}](v0)(v5)

\end{tikzpicture}
\end{center}
\caption{$C_6$}
\end{figure}

When $n$ is even and $U$ is the transition matrix of bipartite walk on $C_n$,
using the same argument as we do when we discuss the transition matrix of bipartite walk on paths, we have that 
\begin{equation}
\label{cycle U permutation}
Ue_i=
\begin{cases}
e_{i+2\Mod n} \quad\text{if }i \text{ is odd;}\\
e_{i-2 \Mod n} \quad\text{if }i \text{ is even.} \\
\end{cases}
\end{equation}

\begin{theorem}
\label{U cycle permutation}
When $n$ is even, the transition matrix $U$ of the bipartite walk on $C_n$ is a cyclic permutation matrix of order $n/2$.
\end{theorem}

\proof The mapping relation~\ref{cycle U permutation} implies that $U$ is a cyclic permutation whose cycle form is  
\[
	(e_0,e_{n-2},\cdots,e_2)(e_1,e_3,\cdots,e_{n-1}).\qed
\]

Note that eigenvalues of $C_n$ are 
\[
	\Bigg\{2\cos\left(\frac{2\pi k}{n}\right): k\in\{0,1,\cdots,n-1\}\Bigg\}.
	\]
So when $n\equiv 0 \Mod 4$, the adjacency matrix of $C_n$ is not invertible and we consider the Hamiltonian of $U^2$ instead.

\begin{corollary}
\label{H-digraph of cycles}
Let $U$ be the transition matrix of of bipartite walk on $C_n$ for some even $n$.
When $n\equiv 2 \Mod 4$, let $H$ be the Hamiltonian of $U$, then the corresponding $H$-digraph is two copies of a weighted oriented $K_{\frac{n}{2}}$. When $n\equiv 0 \Mod 4$ and $n\geq 12$, let $H$ be the Hamiltonian of $U^2$, then the corresponding $H$-digraph is three copies of a weighted oriented $K_{\frac{n}{4}}$.
\end{corollary}

\proof From Theorem~\ref{U cycle permutation}, the transition matrix of $U$ is two $\frac{n}{2}$-cycles and each cycle is the permutation associated with the transition matrix of bipartite walk on $P_{\frac{n}{2}+1}$. Results follow from Theorem~\ref{H-digraph of Path} and Corollary~\ref{odd path H-digraph}.\qed
Note that when $n=4$, the Hamiltonian of $U$ is zero matrix. When $n=8$, the transition matrix $U$ and $U^2$ both have $-1$ as eigenvalues. There is no real skew-symmetric $S$ such that the Hamiltonian of $U$ or the Hamiltonian of $U^2$ is of the form $iS$ and so, we omit the case when $n=8$.

\section{Universal PST}

Let $U$ be the transition matrix of the continuous walk defined over graph $G$, then we say there is perfect state transfer from state $a$ to state $b$ if \[
	\abs{U(t)_{a,b}}=1.
\] 
A graph $G$ has universal perfect state transfer if it has perfect state transfer between every pair of its vertices.
According to Cameron et al.~in~\cite{Cameron2014}, the only known graphs that have universal perfect state transfer are oriented $K_2,C_3$ with constant weight $i$ assigned on each arc.

In this section, we show that bipartite walk can help us to construct weighted oriented graphs where the continuous quantum walk has universal perfect state transfer. Note that when we talk about continuous walks on weighted graph, the Hamiltonian is the weighted adjacency matrix $A$ of the graph, i.e., the transition matrix is of the form 
\[
	\exp(iA).
	\] 

If the transition matrix $U$ of a bipartite walk is a permutation matrix with finite order, then its $H$-digraph has universal perfect state transfer.

\begin{lemma} 
\label{when U is permutation}
Let $G$ be a connected bipartite walk.
The transition matrix of the bipartite walk on $G$ is a permutation matrix if and only if every vertex of $G$ has degree either $1$ or $2$.
\end{lemma}

\proof Here, we use the same notations as defined in Section~\ref{Intro}. If every vertex of $G$ has degree either $1$ or $2$, using the same notations as before, then both $2P-I$ and $2Q-I$ are permutation matrices. Hence, the transition matrix $U$ is also a permutation matrix. 

For the other direction, note that $2P-I,2Q-I$ are reflections about the spaces spanned by characteristic vectors of cells of $\pi_0,\pi_1$ respectively and cells in one partition are disjoint. Then in order for $U$ to map an edge $e_i$ to another edge $e_j$, the size of each cell of both partitions $\pi_1,\pi_2$ cannot be greater than two.\qed

We have shown in Theorem~\ref{path U permutation} that the transition matrix of the bipartite walk over $P_n$ for some even $n$  is a permutation matrix with finite order. We can use this to produce weighted graphs over which continuous walks have universal perfect state transfer.

The following theorem follows directly from the fact that $U^{n-1}=I$ and Theorem~\ref{H-digraph of Path}.

\begin{corollary} Let $n$ be an even integer. Let $s,t$ be distinct integer in $\{0,\cdots,n-2\}$. we define 
\[
\alpha=\begin{cases}
 \frac{t-s}{2},\quad\text{if both }s, t \text{ are odd;}\\[2mm]
\frac{s+t+1}{2} \quad\text{if }s\text{ is even and } t \text{ is odd;}\\[2mm]
\frac{-t-s-1}{2},\quad\text{if }s\text{ is odd and } t \text{ is even;}\\[2mm]
 \frac{s-t}{2},\quad\text{if both }s, t \text{ are even.}\\[2mm]
 \end{cases}.
\] 
When $n$ is even, the edge $(s,t)$ of $K_{n-1}$ is assigned with weight 
\[
\frac{2}{n-1}
\sum_{r=1}^{\frac{n}{2}-1}\frac{2\pi r}{(n-1)} \sin\left(\frac{2\pi r}{n-1}\alpha\right)
\] 
for all distinct $s,t\in\{0,\cdots,n-2\}$. Let $A$ be the weighted adjacency matrix of the resulting weighted $K_{n-1}$. Then the continuous walk with transition matrix $\exp(iA)$ has universal perfect state transfer and every state will get transferred perfectly to any other state within time $t\leq n-1$.
\end{corollary}

\section{Open questions}
Since continuous quantum walks whose Hamiltonians are symmetric, perfect state transfer is symmetric. That is, in continuous walks, there exists time $t$ when there is perfect state transfer from state $a$ to $b$ and from state $b$ to $a$. However, perfect state transfer in the discrete quantum walk is not necessarily symmetric. Because the transition matrices of discrete quantum walks are not symmetric in general, there is no guarantee that there exists a positive integer $k$ such that at $k$-th step there is perfect state transfer between two states. In fact, there may be cases where there is perfect state transfer from state $a$ to state $b$ while there is no perfect state transfer from state $b$ to state $a$.
 
Recall that the transition matrix of the bipartite walk defined on the graph in Figure~\ref{not return pst graph} is 
\[
U=\begin{pmatrix}
0 & -\frac{1}{3} & 0 &  \frac{2}{3}&  \frac{2}{3} & 0 & 0 \\[2.5mm]
0 & 0 & 0 & 0 & 0 & 0 & 1 \\[2.5mm]
0 & \frac{2}{3}  & 0 & \frac{2}{3}  & \frac{1}{3} & 0 & 0 \\[2.5mm]
0 &  \frac{2}{3}  & 0 & -\frac{1}{3}  &  \frac{2}{3}  & 0 & 0 \\[2.5mm]
0 & 0 & 1 & 0 & 0 & 0 & 0 \\[2.5mm]
1 & 0 & 0 & 0 & 0 & 0 & 0 \\[2.5mm]
0 & 0 & 0 & 0 & 0 & 1 & 0
\end{pmatrix}
.\]
State $e_i$ is the characteristic vector of $i$. It is easy to see that there is perfect state transfer from state $e_1$ to $e_6$ at step $k=1$. But up to $k=300000$ steps, there is no perfect state transfer observed from $e_6$ to $e_1$. We suspect that there is no perfect state transfer from $e_6$ to $e_1$. We would like to find a condition on graph $G$ that determines whether or not perfect state transfer is symmetric.

So far, the graphs we observed, over which bipartite walks defined has perfect state transfer, all have minimum degree at most two. We would like to know if there is any graph $G$ with minimum degree at least three that has perfect state transfer in the bipartite walk defined on $G$.

We would like to know how the structure of the graph $G$ affects behaviors of state transfer in the bipartite walk and if there is any feature of bipartite walk that can be determined by the combinatorial or algebraic properties of the graph it is defined on. This will be the future direction of our studies.

\bibliographystyle{plain}

\end{document}